\documentclass[11pt]{article}
\usepackage{pifont}

\usepackage{amsfonts}
\usepackage{latexsym}
\usepackage{amsmath}
\usepackage{amssymb}
\usepackage{color}

\newtheorem{theorem}{Theorem}[section]

\newtheorem{corollary}{Corollary}[section]

\newenvironment{proof}[1][Proof]{\noindent \textbf{#1.} }{\ \ \  $\Box$}

\newtheorem{lemma}{Lemma}[section]
\newtheorem{definition}{Definition}[section]
\newtheorem{remark}{Remark}[section]

\title{Solvability of general backward stochastic Volterra integral
 equations with non-Lipschitz conditions \thanks{This work is supported by National Natural
Science Foundation of China Grant 10771122, Natural Science
Foundation of Shandong Province of China Grant Y2006A08 and
National Basic Research Program of China (973 Program, No.
2007CB814900).}}

\date{January 18 2010}

\author{Tianxiao Wang and Yufeng Shi \thanks{Corresponding author, E-mail:yfshi@sdu.edu.cn}\\ \small{School of
Mathematics, Shandong University, Jinan 250100, China}}

\begin{document}

\maketitle

\begin{abstract}
In this paper we study the unique solvability of backward stochastic
Volterra integral equations (BSVIEs in short), in terms of both the
M-solutions introduced in \cite{Y2} and the adapted solutions in
\cite{L}, \cite{R} or \cite{WZ}. A general existence and uniqueness
of M-solutions is proved under non-Lipschitz conditions by virtue of
a briefer argument than the one in \cite{Y2}, which extends the
results in \cite{Y2}. For the adapted solutions, the unique
solvability of BSVIEs under more general stochastic non-Lipschitz
conditions is obtained, which generalize the results in \cite{L},
\cite{R} and \cite{WZ}.
\par  $\textit{Keywords:}$ Backward stochastic Volterra integral
equations,
   Adapted M-solutions, Non-Lipschitz condition, stochastic Lipschitz coefficients, adapted solutions
\end{abstract}



\section{Introduction}
\label{sec:intro}
Let $\{W_t\}_{t\in [0,T]}$ be a $d$-dimensional Wiener process
defined on a probability space $(\Omega ,\mathcal {F},P)$ and
$\{\mathcal{F}_t\}_{t\in [0,T]} $ denote the natural filtration of
$\{W_t\},$ such that $\mathcal{F}_0$ contains all $P$-null sets of
$\mathcal{F}.$ This paper is motivated by the recent work of Yong
(\cite{Y1}, \cite{Y2}), which studied an extension of backward
stochastic differential equations (BSDEs in short), i.e. backward
stochastic Volterra integral equations (BSVIEs in short). The
nonlinear BSDEs of the form
\begin{equation}
Y(t)=\xi +\int_t^Tg(s,Y(s),Z(s))ds-\int_t^TZ(s)dW(s),
\end{equation}
initiated by Pardoux and Peng \cite{PP1}, have been studied
extensively in the past two decades. We refer the reader to the
books of Ma and Yong \cite{MY}, Yong and Zhou \cite{YZ} and the
survey paper of El Karoui, Peng and Quenez \cite{EPQ} for the
detailed accounts of both theory and application (especially in
mathematical finance and stochastic control) for such equations. On
the other hand, BSVIEs of the form
\begin{equation}
Y(t)=\psi
(t)+\int_t^Tg(t,s,Y(s),Z(t,s),Z(s,t))ds-\int_t^TZ(t,s)dW(s),
\end{equation}
were firstly introduced by Yong \cite{Y1}. We refer the reader to
\cite{Y1}, \cite{Y3} and \cite{Y2} for both theory and application
in dynamic risk measure and optimal control. As to the adapted
solution of BSVIE (2) ($g$ is independent of $Z(s,t)$ or
$\psi(t)=\xi$), see \cite{A}, \cite{L}, \cite{R}, \cite{WZ}, and the
references cited therein.

No matter the M-solution in \cite{Y1}, \cite{Y2} and \cite{Y3}, or
the adapted solution in \cite{L}, \cite{R} and \cite{WZ}, they made
at least one of the following assumptions, 1) $g$ is independent of
$Z(s,t)$, 2) the terminal condition is $\mathcal{F}_{T}$-measurable
random variable $\xi$, 3) the Lipschitz condition, moreover, the
coefficient is deterministic, 4) the deterministic non-Lipschitz
condition. In this paper, we consider the general case of the above
two kinds of solutions respectively. At first we will prove the
unique solvability of M-solutions with a new method. The reason is
at least two-fold. On the one hand, before proving the unique
existence of M-solution, we should make many preparations if we use
the method in \cite{Y2}, such as the solvability of solutions of
certain stochastic Fredholm integral equation and some other
estimates of M-solution for certain simple BSVIEs. On the other
hand, BSVIEs do not have time-consistency (or semigroup) property,
and the process $Z$ has two parameters, so if we use induction
method in the solution procedure for BSVIEs, we have to use four
steps as in \cite{Y2}, which seems rather complicated and
sophisticated. So we will introduce a new convenient method from
other perspective, that is, when
$T$ is finite, we can use an equivalent norm in $\mathcal{H}%
^2[0,T]$ as follows:
\[
\left\| (y(\cdot ),z(\cdot ,\cdot ))\right\|
_{\mathcal{H}^2[0,T]}=\left[ E\int_0^Te^{\beta
t}|y(t)|^2dt+E\int_0^T\int_0^Te^{\beta s}|z(t,s)|^2dsdt\right]
^{\frac 12},
\]
where $\beta$ is a positive constant, and then we will prove the
results of M-solutions within the new norm with one step. This is
our first contribution in this paper. When $g$ is independent of
$Z(s,t),$ we can use the estimate in lemma 3.1 and the method in
\cite{WZ} to get the solvability of M-solution of BSVIE (2) ($g$ is
independent of $Z(s,t)$) under non-Lipschitz condition. However, as
to the general form of BSVIE (2), the method in \cite{WZ} does not
work any more because of the appearance of $E\int_t^T|Z(s,t)|^2ds$,
which can be estimated by means of Malliavin calculus, see
\cite{Y2}, and this will make the problem more complicated. So we
have to replace $e^{\beta t}E|Y(t)|^2+E\int_t^Te^{\beta
s}|Z(t,s)|^2ds$, as in \cite{WZ}, with a weaker form
$\int_u^Te^{\beta t}E|Y(t)|^2dt+E\int_u^T\int_t^Te^{\beta
s}|Z(t,s)|^2dsdt$, $u\in[0,T].$ We also have to prove a new lemma by
means of the definition of concave function, then we can obtain the
unique existence of M-solution of BSVIE (2) under non-Lipschitz
condition, which generalize the result in \cite{Y1}, \cite{Y2} and
\cite{Y3}. Thirdly, we claim that It$\hat{o}$ formula plays a key
role in the BSDEs case, as well as the BSVIEs case in \cite{WZ}. One
question is can we get the solvability of adapted solution of (2)
($g$ is independent of $Z(s,t)$) under stochastic non-Lipschitz
conditions without involving It$\hat{o}$ formula? The answer is
positive and we will prove it in the following, which generalize the
result in \cite{L} and \cite{WZ}.

Recently the author considered the unique solvability of M-solution
under non-Lipschitz condition by induction in \cite{R}. Note that
our method here is different from it, moreover, briefer than it. On
the other hand, the assumption on the coefficients in \cite{R} is
also much stronger than ours here.

The paper is organized as follows. In Section 2, we give some
preliminary results and notations which are needed in the following
sections. An important estimate for M-solutions (or adapted
solutions) is presented in Subsection 3.1. With this estimate, we
give the existence and uniqueness result of M-solutions under
Lipschitz condition in Subsection 3.2. The case of adapted solutions
is also treated. In Subsection 3.3, we consider the unique
solvability of M-solutions (adapted solutions respectively) under
non-Lipschitz case. At last examples of M-solutions and adapted
solutions under non-Lipschitz condition is present.

\section{Preliminaries}

   In this section, we will make some preliminaries. $\forall R,S\in[0,T],$ in the
following we denote $\Delta ^c[R,S]=\{(t,s)\in[R,S]^{2}; t\leq s\},$
$\Delta^c=\Delta^c[0,T],$ $\Delta [R,S]=\{(t,s)\in[R,S]^{2}; t>s\},$
and $\Delta=\Delta[0,T].$ Let us first introduce some spaces. Let
$\beta $ be a positive constant. $ A(t)$ is a non-negative $\mathcal
{F}_{t}$-adapted increasing process. Let $L_{\mathcal {F}_{T}}
^{2,\beta}[0,T]$ be the set of the
$\mathcal{B}([0,T])\otimes\mathcal{F}_{T}$-measurable processes
$X:[0,T]\times \Omega \rightarrow R^m$ such that $ Ee^{\beta
A(T)}\displaystyle\int_0^T|X(t)|^2dt<\infty. $ We denote
\begin{eqnarray*}
\mathcal{H} ^{2,\beta}[R,S] &=&L_{\mathcal {F}}^{2,\beta}[R,S]\times
L
^{2,\beta}(R,S;L_{\mathcal {F}}^2[R,S]), \\
\mathcal{H}_{t}^{2,\beta}[R,S] &=&L_{\mathcal {F}}^{2,\beta
}[R,S]\times L ^{2,\beta}(R,S;L_{\mathcal {F}}^2[t,S]).
\end{eqnarray*}
Here $L_{\mathcal {F}}^{2,\beta}[R,S]$ is the set of all adapted
processes $X:[R,S]\times \Omega \rightarrow R^m$ such that $
E\displaystyle\int_R^Se^{\beta A(s)}|X(s)|^2ds<\infty . $
$L ^{2,\beta}(R,S;L_{\mathcal {F}}^2[R,S])$ is the set of all processes $%
Z:[R,S]^2\times \Omega \rightarrow R^{m\times d}$ such that for almost all $%
t\in [R,S],$ $Z(t,\cdot )$ is $\mathcal {F}$-progressively
measurable satisfying $
E\displaystyle\int_R^S\displaystyle\int_R^Se^{\beta
A(s)}|Z(t,s)|^2dsdt<\infty . $
$L ^{2,\beta}(R,S;L_{\mathcal {F}}^2[t,S])$ is the set of all processes $%
Z(t,s):\Delta^c[R,S]\times \Omega \rightarrow R^{m\times d}$ such
that for almost all $t\in [R,S],$ $Z(t,\cdot )$ is $\mathcal
{F}$-progressively measurable satisfying $
E\displaystyle\int_R^S\displaystyle\int_t^Se^{\beta
A(s)}|Z(t,s)|^2dsdt<\infty . $ Let $L_{\mathcal {F}_{T}}^2[0,T]$ be
the set of the $\mathcal{B}([0,T])\otimes\mathcal{F}_{T}$ processes
$X:[0,T]\times \Omega \rightarrow R^m$ such that $
E\displaystyle\int_0^T|X(t)|^2dt<\infty. $ We also denote
\begin{eqnarray*}
\mathcal{H}^2[R,S]=L_{\mathcal {F}}^2[R,S]\times L^2(R,S;L_{\mathcal
{F}}^2[R,S]), \\
\mathcal{H}_{t}^2[R,S] = L_{\mathcal {F}}^2[R,S]\times L
^2(R,S;L_{\mathcal {F}}^2[t,S]).
\end{eqnarray*}
Here $L_{\mathcal {F}}^2[R,S]$ is the set of all adapted processes
$X:[R,S]\times \Omega \rightarrow R^m$ such that $
E\displaystyle\int_R^S|X(s)|^2ds<\infty . $
$L^2(R,S;L_{\mathcal {F}}^2[R,S])$ is the set of all processes $%
Z:[R,S]^2\times \Omega \rightarrow R^{m\times d}$ such that for almost all $%
t\in [R,S],$ $Z(t,\cdot )$ is $\mathcal {F}$-progressively
measurable satisfying $
E\displaystyle\int_R^S\int_R^S|Z(t,s)|^2dsdt<\infty . $
$L^2(R,S;L_{\mathcal {F}}^2[t,S])$ is the set of all processes $%
Z(t,s):\Delta^c[R,S]\times \Omega \rightarrow R^{m\times d}$ such
that for almost all $t\in [R,S],$ $Z(t,\cdot )$ is $\mathcal
{F}$-progressively measurable satisfying $
E\displaystyle\int_R^S\int_t^S|Z(t,s)|^2dsdt<\infty . $ Now we give
two definitions needed in the sequel.
\begin{definition}
Let $S\in [0,T]$. A pair of $(Y(\cdot ),Z(\cdot ,\cdot ))\in \mathcal{H}%
 ^{2,\beta}[S,T]$ is called an adapted $M$-solution of BSVIE (2) on
$[S,T]$ if (2) holds in the usual It\^o's sense for almost all $t\in
[S,T]$ and, in addition, the following holds:
\[
Y(t)=E^{\mathcal{F}_S}Y(t)+\int_S^tZ(t,s)dW(s),\quad t\in [S,T].
\]
\end{definition}
\begin{definition}
A pair of $(Y(\cdot ),Z(\cdot ,\cdot ))\in \mathcal
{H}_{t}^{2,\beta}[0,T]$ is called an adapted solution of the
following simple BSVIE (3) if (3) holds in the usual It\^o's sense
\begin{equation}
Y(t)=\psi (t)+\int_t^Tg(t,s,Y(s),Z(t,s))ds-\int_t^TZ(t,s)dW(s),\quad
t\in [0,T].
\end{equation}
\end{definition}
In \cite{Y2}, the author gave the definition of M-solution of BSVIE in $\mathcal{H%
}^2[0,T]$. The author also considered the existence and uniqueness
of adapted solution of (3) ($\psi(\cdot)$ is
replaced with $\xi$) in $\mathcal{H%
}^2_t[0,T]$ in \cite{WZ}.

We give the following assumptions of $g$ for BSVIE (2):

(H1) Let $g:\Delta ^c\times R^m\times R^{m\times d}\times R^{m\times
d}\times \Omega \rightarrow R^m$ be $\mathcal{B}(\Delta ^c\times
R^m\times R^{m\times
d}\times R^{m\times d})\otimes \mathcal{F}_T$-measurable such that $%
s\rightarrow g(t,s,y,z,\zeta )$ is $\mathcal {F}$-progressively
measurable for all $(t,y,z,\zeta )\in [0,T]\times R^m\times
R^{m\times d}\times R^{m\times d}$, furthermore, $g$ satisfies the
Lipschitz conditions with stochastic coefficient, i.e., $\forall y,$
$\overline{y}\in R^m,$ $z,$ $\overline{z},$ $\zeta ,$
$\overline{\zeta }\in R^{m\times d},$
\begin{eqnarray*}
&&|g(t,s,y,z,\zeta )-g(t,s,\overline{y},\overline{z},\overline{\zeta })| \\
&\leq &L(t,s)(r_1(s)|y-\overline{y}|+r_2(s)|z-\overline{z}|+r_3(s)|\zeta -%
\overline{\zeta }|),
\end{eqnarray*}
where $(t,s)\in \Delta ^c,$  $r_1(s),$ $r_2(s)$ and $r_3(s)$ are
non-negative adapted processes and we denote
\[
\alpha ^2(s)=r_1^2(s)+r_2^2(s)+r_3^2(s),\quad A(t)=\int_0^t\alpha
^2(s)ds.
\]
We assume that $\alpha ^2(s)\geq \delta,$ where $\delta$ is a
positive constant, $\alpha(s)$ is a positive adapted process and
$L(t,s)$ is a deterministic non-negative function. Furthermore, we
assume
\[
E\int_0^T\int_t^Te^{\beta A(s)}|g_0(t,s)|^2dsdt<\infty ,
\]
where $g_0(t,s)=g(t,s,0,0,0).$

\section{Main results for M-solutions}

\subsection{A basic estimate for M-solutions of BSVIEs}

In this subsection, inspired by the method of estimating the adapted
solutions of BSDEs in \cite{EH}, we give a lemma which is needed in
the following.

\begin{lemma}
We consider the following simple BSVIE
\begin{equation}
Y(t)=\psi (t)+\int_t^Tf(t,s)ds-\int_t^TZ(t,s)dW(s),\quad t\in [0,T],
\end{equation}
where $\psi (\cdot )\in L_{\mathcal {F}_{T}}^{2,\beta} [0,T],$
$f:\Delta^c\times\Omega\rightarrow R^m$ be
$\mathcal{B}(\Delta^c)\otimes\mathcal{F}_{T}$-measurable such that
$s\rightarrow f(t,s)$ is $\mathcal{F}$-progressively measurable for
all $t\in[0,T],$ and $E\int_0^T\int_t^Te^{\beta
A(s)}|f(t,s)|^2dsdt<\infty .$ Then (4) admits a unique adapted
solution $(Y(\cdot ),Z(\cdot ,\cdot ))\in \mathcal{H}_{t
}^{2,\beta}[0,T],$ and we have the following estimate:
\begin{eqnarray}
&&\ E\int_0^Te^{\beta A(s)}|Y(s)|^2ds+E\int_0^T\int_t^Te^{\beta
A(s)}|Z(t,s)|^2dsdt  \nonumber \\
\ &\leq &CE\int_0^Te^{\beta A(t)}|\psi (t)|^2dt+CE\int_0^Te^{\beta
A(t)}\left| \int_t^Tf(t,s)ds\right| ^2dt  \nonumber \\
&&+CE\int_0^T\int_t^Te^{\beta A(u)}|\psi (t)|^2d\beta A(u)dt  \nonumber \\
&&+\ CE\int_0^T\int_t^Te^{\beta A(s)}\left| \int_s^Tf(t,u)du\right|
^2d\beta A(s)dt.
\end{eqnarray}
Furthermore,
\begin{eqnarray}
&&\ \ E\int_0^Te^{\beta A(s)}|Y(s)|^2ds+E\int_0^T\int_t^Te^{\beta
A(s)}|Z(t,s)|^2dsdt \nonumber  \\
\ &\leq &CEe^{\beta A(T)}\int_0^T|\psi (t)|^2dt+\frac C\beta
E\int_0^T\int_t^Te^{\beta A(s)}\frac{|f(t,s)|^2}{\alpha ^2(s)}dsdt.
\end{eqnarray}
Hereafter $C$ is a generic positive constant which may be different
from line to line.
\end{lemma}
\begin{proof}We consider a family of BSDEs with parameters $t$ on $[0,T]$ in the
following form:
\begin{equation}
\lambda (t,r)=\psi (t)+\int_r^Tf(t,s)ds-\int_r^T\mu (t,s)dW(s),\quad
t,r\in [0,T].
\end{equation}
By the classical existence and uniqueness theorem of BSDE in
\cite{PP1}, there
exists a unique solution $(\lambda (t,\cdot ),\mu (t,\cdot ))$ for every $%
t\in [0,T].$ Let $Y(t)=\lambda (t,t),\quad Z(t,s)=\mu (t,s),$ $t\leq
s.$ Then we obtain the existence and uniqueness of the adapted
solution for (4). From (7) we arrive at, $\forall r\geq t,$%
\begin{eqnarray}
\lambda (t,r) &=&E^{\mathcal{F}_r}\left( \psi (t)+\int_r^Tf(t,s)ds %
 \right),  \nonumber
 \end{eqnarray}
 and
 \begin{eqnarray}
\int_r^TZ(t,s)dW(s) &=&\int_r^T\mu (t,s)dW(s)=\psi
(t)+\int_r^Tf(t,s)ds-\lambda (t,r).
\end{eqnarray}
Especially when $r=t,$%
\[
\int_t^TZ(t,s)dW(s)=\int_t^T\mu (t,s)dW(s)=\psi
(t)+\int_t^Tf(t,s)ds-Y(t).
\]
Now we estimate $E\displaystyle\int_0^Te^{\beta
A(s)}|Y(s)|^2ds+E\displaystyle\int_0^T\displaystyle\int_t^Te^{\beta
A(s)}|Z(t,s)|^2dsdt.$ By Cauchy-Schwarz inequality we deduce that
\begin{eqnarray}
\left| \int_s^Tf(t,u)du\right| ^2 &=&\left|
\int_s^Te^{\frac{-rA(u)}2}\alpha
(u)e^{\frac{rA(u)}2}\frac{f(t,u)}{\alpha (u)}du\right| ^2   \nonumber  \\
&\leq &\int_s^Te^{-rA(u)}\alpha ^2(u)du\cdot \int_s^Te^{rA(u)}\frac{%
|f(t,u)|^2}{\alpha ^2(u)}du  \nonumber \\
&\leq &\frac 1re^{-rA(s)}\int_s^Te^{rA(u)}\frac{|f(t,u)|^2}{\alpha ^2(u)}%
du,\quad t,s\in [0,T],
\end{eqnarray}
where $r$ is a positive constant. By taking $r=\frac \beta 2$ in
(9), we see that
\begin{eqnarray*}
&&\int_t^Te^{\beta A(s)}\left| \int_s^Tf(t,u)du\right| ^2d\beta A(s) \\
&\leq &\frac 4\beta \int_t^Te^{\frac \beta 2A(s)}\left(
\int_s^Te^{\frac
\beta 2A(u)}\frac{|f(t,u)|^2}{\alpha ^2(u)}du\right) d\frac \beta 2A(s) \\
&=&\left. \frac 4\beta e^{\frac \beta 2A(s)}\left( \int_s^Te^{\frac
\beta
2A(u)}\frac{|f(t,u)|^2}{\alpha ^2(u)}du\right) \right| _t^T \\
&&+\frac 4\beta \int_t^Te^{\beta A(s)}\frac{|f(t,s)|^2}{\alpha ^2(s)}ds \\
&\leq &\frac 4\beta \int_t^Te^{\beta A(s)}\frac{|f(t,s)|^2}{\alpha
^2(s)}ds.
\end{eqnarray*}
Therefore,
\begin{equation}
E\int_0^T\int_t^Te^{\beta A(s)}\left| \int_s^Tf(t,u)du\right|
^2d\beta
A(s)dt\leq \frac 4\beta E\int_0^T\int_t^Te^{\beta A(s)}\frac{|f(t,s)|^2}{%
\alpha ^2(s)}dsdt.
\end{equation}
We also obtain the following result by taking $s=t$ and $r=\beta $
in (9),
\[
E\int_0^Te^{\beta A(t)}\left| \int_t^Tf(t,u)du\right| ^2dt\leq \frac
1\beta E\int_0^T\int_t^Te^{\beta A(u)}\frac{|f(t,u)|^2}{\alpha
^2(u)}dudt.
\]
At first we estimate
$E\displaystyle\int_0^T\displaystyle\int_t^Te^{\beta
A(s)}|Z(t,s)|^2dsdt.$ Obviously, we have $t,r\in [0,T],$
\begin{eqnarray}
&&\int_r^Te^{\beta A(s)}\left( \int_s^T|Z(t,u)|^2du\right) d\beta A(s) \nonumber \\
&=&\left. e^{\beta A(s)}\left( \int_s^T|Z(t,u)|^2du\right) \right|
_r^T+\int_r^Te^{\beta A(s)}|Z(t,s)|^2ds.
\end{eqnarray}
For arbitrary $t\in [0,T],$ we can rewrite (11) after taking $r=t$,
\begin{eqnarray}
&&E\int_0^T\int_t^Te^{\beta A(s)}|Z(t,s)|^2dsdt  \nonumber \\
&=&E\int_0^T\int_t^Te^{\beta A(s)}\left( \int_s^T|Z(t,u)|^2du\right)
d\beta A(s)dt \nonumber \\
&&+E\int_0^Te^{\beta A(t)}\int_t^T|Z(t,u)|^2dudt.
\end{eqnarray}
Now we give a estimate to the second expression in the right part of
(12)
\begin{eqnarray}
&&E\int_0^Te^{\beta A(t)}\int_t^T|Z(t,u)|^2dudt \nonumber \\
&=&E\int_0^TE\left( e^{\beta A(t)}\left. \int_t^T|Z(t,u)|^2du\right|
\mathcal{F}_t\right) dt  \nonumber \\
&=&E\int_0^Te^{\beta A(t)}E\left( \left. \left( \int_t^TZ(t,u)d
W(u)\right)
^2\right| \mathcal{F}_t\right) dt  \nonumber \\
&=&E\int_0^Te^{\beta A(t)}E\left( \left. \left( \psi
(t)+\int_t^Tf(t,u)du-Y(t)\right) ^2\right| \mathcal{F}_t\right) dt
\nonumber
\\
&\leq &3E\int_0^Te^{\beta A(t)}|\psi (t)|^2dt+3E\int_0^Te^{\beta
A(t)}\left|
\int_t^Tf(t,u)du\right| ^2dt  \nonumber \\
&&+3E\int_0^Te^{\beta A(t)}|Y(t)|^2dt  \nonumber \\
&\leq &3E\int_0^Te^{\beta A(t)}|\psi (t)|^2dt+\frac 3\beta
E\int_0^T\int_t^Te^{\beta A(u)}\frac{|f(t,u)|^2}{\alpha ^2(u)}dudt
\nonumber
\\
&&+3E\int_0^Te^{\beta A(t)}|Y(t)|^2dt.
\end{eqnarray}
Obviously, we can use the similar method as (13) to estimate the
first expression in the right part of (12) as follows:
\begin{eqnarray}
&&E\int_0^T\int_t^Te^{\beta A(s)}\left( \int_s^T|Z(t,u)|^2du\right)
d\beta
A(s)dt  \nonumber \\
&\leq &3E\int_0^T\int_t^Te^{\beta A(s)}|\psi (t)|^2d\beta A(s)dt
\nonumber
\\
&&+3E\int_0^T\int_t^Te^{\beta A(s)}\left( \int_s^Tf(t,u)du\right)
^2d\beta
A(s)dt  \nonumber \\
&&+3E\int_0^T\int_t^Te^{\beta A(s)}|\lambda (t,s)|^2d\beta A(s)dt.
\end{eqnarray}
For the second expression in the right part of (14), inequality (10)
implies that
\begin{eqnarray}
&&3E\int_0^T\int_t^Te^{\beta A(s)}\left( \int_s^Tf(t,u)du\right)
^2d\beta
A(s)dt   \nonumber \\
&\leq &\frac{12}\beta E\int_0^T\int_t^Te^{\beta A(s)}\frac{|f(t,s)|^2}{%
\alpha ^2(s)}dsdt.
\end{eqnarray}
It is time for us to estimate the third expression of (14). Observe
that
\[
\lambda (t,s)=E^{\mathcal{F}_s}\left( \psi (t)+\int_s^Tf(t,u)du
\right),\quad t\leq s,
\]
so we deduce that,
\[
|\lambda (t,s)|^2\leq 2E(|\psi (t)|^2\left| \mathcal{F}_s\right.
)+2E\left( \left. \left| \int_s^Tf(t,u)du\right| ^2\right|
\mathcal{F}_s\right) ,
\]
and
\begin{eqnarray}
&&3E\int_0^T\int_t^Te^{\beta A(s)}|\lambda (t,s)|^2d\beta A(s)dt
\nonumber
\\
&\leq &6E\int_0^T\int_t^Te^{\beta A(s)}|\psi (t)|^2d\beta A(s)dt
\nonumber
\\
&&+6E\int_0^T\int_t^Te^{\beta A(s)}\left| \int_s^Tf(t,u)du\right|
^2d\beta
A(s)dt  \nonumber \\
&\leq &6E\int_0^T|\psi (t)|^2dt\left( \int_t^Te^{\beta A(s)}d\beta
A(s)\right)  \nonumber \\
&&+\frac{24}\beta E\int_0^T\int_t^Te^{\beta
A(s)}\frac{|f(t,s)|^2}{\alpha ^2(s)}dsdt.
\end{eqnarray}
Equality (14), (15) and (16) implies that
\begin{eqnarray}
&&E\int_0^T\int_t^Te^{\beta A(s)}\left( \int_s^T|Z(t,u)|^2du\right)
d\beta
A(s)dt   \nonumber \\
&\leq &\frac{36}\beta E\int_0^T\int_t^Te^{\beta A(s)}\frac{|f(t,s)|^2}{%
\alpha ^2(s)}dsdt  \nonumber \\
&&+9E\int_0^T|\psi (t)|^2dt\left( \int_t^Te^{\beta A(s)}d\beta
A(s)\right) .
\end{eqnarray}
From (12), (13) and (17), we also see that
\begin{eqnarray*}
&&E\int_0^T\int_t^Te^{\beta A(s)}|Z(t,s)|^2dsdt \\
&\leq &3E\int_0^Te^{\beta A(t)}|\psi (t)|^2dt+9E\int_0^T|\psi
(t)|^2dt\left(
\int_t^Te^{\beta A(s)}d\beta A(s)\right) \\
&&+\frac{39}\beta E\int_0^T\int_t^Te^{\beta
A(s)}\frac{|f(t,s)|^2}{\alpha ^2(s)}dsdt+3E\int_0^Te^{\beta
A(t)}|Y(t)|^2dt.
\end{eqnarray*}
Due to $Y(t)=E^{\mathcal{F}_t}\left( \psi
(t)+\displaystyle\int_t^Tf(t,s)ds \right) $, we have
\begin{eqnarray}
&&E\int_0^Te^{\beta A(t)}|Y(t)|^2dt \nonumber  \\
&\leq &2E\int_0^Te^{\beta A(t)}|\psi (t)|^2dt+2E\int_0^Te^{\beta
A(t)}\left|
\int_t^Tf(t,s)ds\right| ^2dt  \nonumber \\
&\leq &2E\int_0^Te^{\beta A(t)}|\psi (t)|^2dt+\frac 2\beta
E\int_0^T\int_t^Te^{\beta A(s)}\frac{|f(t,s)|^2}{\alpha ^2(s)}dsdt.
\end{eqnarray}
Eventually we obtain:
\begin{eqnarray*}
&&E\int_0^Te^{\beta A(s)}|Y(s)|^2ds+E\int_0^T\int_t^Te^{\beta
A(s)}|Z(t,s)|^2dsdt \\
&\leq &11E\int_0^Te^{\beta A(t)}|\psi (t)|^2dt+9E\int_0^T|\psi
(t)|^2dt\left( \int_t^Te^{\beta A(u)}d\beta A(u)\right) \\
&&+11E\int_0^Te^{\beta A(t)}\left| \int_t^Tf(t,s)ds\right|
^2dt+9E\int_0^T\int_t^Te^{\beta A(s)}\left| \int_s^Tf(t,u)du\right|
^2d\beta A(s)dt.
\end{eqnarray*}
Furthermore, it follows that:
\begin{eqnarray*}
&&\ E\int_0^Te^{\beta A(t)}|Y(t)|^2dt+E\int_0^T\int_t^Te^{\beta
A(s)}|Z(t,s)|^2dsdt \\
&\leq &20Ee^{\beta A(T)}\int_0^T|\psi (t)|^2dt+\frac{47}\beta
E\int_0^T\int_t^Te^{\beta A(s)}\frac{|f(t,s)|^2}{\alpha ^2(s)}dsdt.
\end{eqnarray*}
\end{proof}
\begin{remark}
If we define $Z(t,s), (0\leq s<t\leq T)$ by the following relation
$$Y(t)=E^{\mathcal{F}_{S}}Y(t)+\int_S^tZ(t,s)dW(s), \quad t\in[S,T], \quad \forall
S\in[0,T].$$ Then BSVIE (4) admits a unique M-solution in
$\mathcal{H}^{2,\beta}[0,T].$
\end{remark}
\subsection{The Lipschitz case}

In this subsection, we give the existence and uniqueness of
M-solution under Lipschitz condition with a much more convenient
method. We require $r_i(s)$ to be deterministic functions
$(i=1,2,3)$. First we give a theorem when $L(t,s)$ is bounded.
\begin{theorem}
Let (H1) hold, $\psi(\cdot)\in L^{2,\beta}_{\mathcal{F}_{T}}[0,T]$
and $r_i(s)$ $(i=1,2,3)$ are deterministic functions, $L(t,s)$ is
bounded, then (2) admits a unique M-solution in $\mathcal{H}%
 ^{2,\beta}[0,T]$.
\end{theorem}
\begin{proof}When $A(\cdot)$ is deterministic function, by the
definition of M-solution we see that
\begin{eqnarray}
&&\ E\int_0^T\int_0^te^{\beta A(s)}|Z(t,s)|^2dsdt  \nonumber \\
\ &\leq &E\int_0^Te^{\beta A(t)}\int_0^t|Z(t,s)|^2dsdt  \nonumber \\
\ &=&\int_0^Te^{\beta A(t)}E\int_0^t|Z(t,s)|^2dsdt \nonumber \\
\ &\leq &E\int_0^Te^{\beta A(t)}|Y(t)|^2dt.
\end{eqnarray}
Then
\begin{eqnarray}
&&E\int_0^Te^{\beta A(t)}|Y(t)|^2dt+E\int_0^T\int_0^Te^{\beta
A(s)}|Z(t,s)|^2dsdt   \nonumber \\
&\leq &2E\int_0^Te^{\beta A(t)}|Y(t)|^2dt+E\int_0^T\int_t^Te^{\beta
A(s)}|Z(t,s)|^2dsdt  \nonumber \\
&\leq &CE\int_0^Te^{\beta A(t)}|\psi (t)|^2dt+CE\int_0^Te^{\beta
A(t)}\left|
\int_t^Tf(t,s)ds\right| ^2dt  \nonumber \\
&&+CE\int_0^T|\psi (t)|^2\left( \int_t^Te^{\beta A(u)}d\beta
A(u)\right) dt
\nonumber \\
&&+CE\int_0^T\int_t^Te^{\beta A(s)}\left| \int_s^Tf(t,u)du\right|
^2d\beta
A(s)dt  \nonumber \\
&\leq &CEe^{\beta A(T)}\int_0^T|\psi (t)|^2dt+\frac C\beta
E\int_0^T\int_t^Te^{\beta A(s)}\frac{|f(t,s)|^2}{\alpha ^2(s)}dsdt.
\end{eqnarray}
Let $\mathcal{M} ^{2,\beta}[0,T]$ be the space of all $(y(\cdot
),z(\cdot ,\cdot ))\in \mathcal{H} ^{2,\beta}[0,T]$ such that
\[
y(t)=Ey(t)+\int_0^tz(t,s)dW(s), \quad t\in [0,T].
\]
Clearly, it is a nonempty closed subspace of $\mathcal{H}
^{2,\beta}[0,T].$ Now we consider the following BSVIE:
\begin{equation}
Y(t)=\psi
(t)+\int_t^Tg(t,s,y(s),z(t,s),z(s,t))ds-\int_t^TZ(t,s)dW(s),\quad
t\in [0,T],
\end{equation}
for any $\psi (\cdot )\in L_{T }^{2,\beta}[0,T]$ and $(y(\cdot
),z(\cdot ,\cdot ))\in \mathcal{M} ^{2,\beta}[0,T].$ Hence by Remark
3.1 we know that (21) admits a unique M-solution $(Y(\cdot ),Z(\cdot
,\cdot
))\in \mathcal{M} ^{2,\beta}[0,T],$ and we can define a map $\Theta :\mathcal{M%
} ^{2,\beta}[0,T]\rightarrow \mathcal{M} ^{2,\beta}[0,T]$ by
\[
\Theta (y(\cdot ),z(\cdot ,\cdot ))=(Y(\cdot ),Z(\cdot ,\cdot
)),\quad \forall (y(\cdot ),z(\cdot ,\cdot ))\in \mathcal{M}
^{2,\beta}[0,T].
\]
Let $(\overline{y}(\cdot ),\overline{z}(\cdot ,\cdot ))\in
\mathcal{M}
^{2,\beta}[0,T]$ and $\Theta (\overline{y}(\cdot ),\overline{z}(\cdot ,\cdot ))=(%
\overline{Y}(\cdot ),\overline{Z}(\cdot ,\cdot )).$ From (20) we
deduce that,
\begin{eqnarray*}
&&\ \ E\int_0^Te^{\beta A(t)}|Y(t)-\overline{Y}(t)|^2dt+E\int_0^T\int_0^Te^{%
\beta A(s)}|Z(t,s)-\overline{Z}(t,s)|^2dsdt \\
\ &\leq &CE\int_0^Te^{\beta A(t)}\left| \int_t^TH(t,s)ds\right|
^2dt+CE\int_0^T\int_t^Te^{\beta A(s)}\left| \int_s^TH(t,u)du\right|
^2d\beta
A(s)dt \\
\ &\leq &CE\int_0^Te^{\beta A(t)}\left| \int_t^TG(t,s)ds\right|
^2dt+CE\int_0^T\int_t^Te^{\beta A(s)}\left| \int_s^TG(t,u)du\right|
^2d\beta
A(s)dt \\
\ &\leq &\frac C\beta E\int_0^T\int_t^Te^{\beta A(s)}L^2(t,s)|y(s)-\overline{%
y}(s)|^2dsdt \\
&&+\frac C\beta E\int_0^T\int_t^Te^{\beta A(s)}L^2(t,s)|z(t,s)-\overline{z}%
(t,s)|^2dsdt \\
&&+\frac C\beta E\int_0^T\int_t^Te^{\beta A(s)}L^2(t,s)|z(s,t)-\overline{z}%
(s,t)|^2dsdt \\
\ &\leq &\frac C\beta E\int_0^Te^{\beta A(s)}|y(s)-\overline{y}%
(s)|^2ds+\frac C\beta E\int_0^T\int_t^Te^{\beta A(s)}|z(t,s)-\overline{z}%
(t,s)|^2dsdt,
\end{eqnarray*}
where
\begin{eqnarray*}
H(t,s) &=&g(t,s,y(s),z(t,s),z(s,t))-g(t,s,\overline{y}(s),\overline{z}(t,s),%
\overline{z}(s,t)), \\
G(t,s) &=&L(t,s)(r_1(s)|y(s)-\overline{y}(s)| \\
&&+r_2(s)|z(t,s)-\overline{z}(t,s)|+r_3(s)|z(s,t)-\overline{z}(s,t)|).
\end{eqnarray*}
The fourth inequality above holds from that:
\begin{eqnarray}
&&E\int_0^T\int_t^Te^{\beta A(s)}|z(s,t)-\overline{z}(s,t)|^2dsdt   \nonumber \\
&=&E\int_0^Te^{\beta
A(t)}\int_0^t|z(t,s)-\overline{z}(t,s)|^2dsdt\leq E\int_0^Te^{\beta
A(t)}|y(t)-\overline{y}(t)|^2dt.
\end{eqnarray}
Choosing a sufficient large number $\beta $ so that $\frac C\beta
<1,$ then the mapping $\Theta$ is contracted from $\mathcal{H}
^{2,\beta}[0,T]$ onto itself. Thus there exists a unique fixed point
which is the unique M-solution of BSVIE (2).
\end{proof}
\begin{remark}
From expression (19) and (22), we can know the reason for assuming
$r_i(s)$ to be deterministic. When $A(s)$ is bounded (or continuous
and deterministic) the norm of $\mathcal{H} ^{2,\beta}[0,T]$ is
equivalent to the norm of $\mathcal{H}^2[0,T],$ then BSVIE (2)
admits a unique M-solution in $\mathcal{H}^2[0,T].$ As a result, we
can finish the proof of uniqueness and existence of M-solution in
$\mathcal{H}^2[0,T]$ with only one step, which is much more
convenient than the four steps in \cite{Y2}.
\end{remark}
In Theorem 3.1 the assumption on the coefficient $L(t,s)$ is
stronger than the one in \cite{Y2}. Next we will see that we can
relax it, but we need introduce a new non-negative process
$A^{*}(t)$ and another assumption on $\alpha ^2(t)$. We denote
$A^{*}(t)=\int_0^t\alpha ^{\frac{2p}{2-p}}(s)ds$ $(1<p<2),$
furthermore, we assume that $\alpha ^2(s)\geq1, t\in[0,T]$.
Obviously we have $A^{*}(t)\geq A(t), t\in[0,T].$ We will obtain
some estimates which are important in the proof of the following
theorem. For arbitrary $1<p<2,$ we obtain
\begin{eqnarray}
&&\ \ \left( \int_r^T|f(t,s)|^pds\right) ^{\frac 2p}  \nonumber \\
\ &=&\left( \int_r^Te^{-\tau A^{*}(s)}\alpha ^p(s)e^{\tau A^{*}(s)}\frac{%
|f(t,s)|^p}{\alpha ^p(s)}ds\right) ^{\frac 2p}  \nonumber \\
\ &\leq &\left( \int_r^Te^{-\tau A^{*}(s)\frac 2{2-p}}\alpha ^{\frac{2p}{2-p}%
}(s)ds\right) ^{\frac{2-p}p}\left( \int_r^Te^{\tau A^{*}(s)\frac 2p}\frac{%
|f(t,s)|^2}{\alpha ^2(s)}ds\right)  \nonumber \\
\ &\leq &\left( \frac{2-p}{2\tau }\right) ^{\frac{2-p}p}e^{-\tau
A^{*}(r)\frac 2p}\int_r^Te^{\tau A^{*}(s)\frac
2p}\frac{|f(t,s)|^2}{\alpha ^2(s)}ds, \quad t,r\in [0,T],
\end{eqnarray}
and
\begin{eqnarray}
\left| \int_s^Tf(t,u)du\right| ^2 &=&\left| \int_s^Te^{\frac{-rA^{*}(u)}%
2}\alpha (u)e^{\frac{rA^{*}(u)}2}\frac{|f(t,u)|}{\alpha
(u)}du\right| ^2 \nonumber
\\
\ &\leq &\int_s^Te^{-rA^{*}(u)}\alpha ^2(u)du\cdot \int_s^Te^{rA^{*}(u)}%
\frac{|f(t,u)|^2}{\alpha ^2(u)}du  \nonumber \\
\ &\leq &\int_s^Te^{-rA^{*}(u)}\alpha ^{\frac{2p}{2-p}}(u)du\cdot
\int_s^Te^{rA^{*}(u)}\frac{|f(t,u)|^2}{\alpha ^2(u)}du  \nonumber \\
\ &\leq &\frac
1re^{-rA^{*}(s)}\int_s^Te^{rA^{*}(u)}\frac{|f(t,u)|^2}{\alpha
^2(u)}du,\quad t,s\in [0,T].
\end{eqnarray}
In (23), let $\tau =\frac p2\beta ,$ we arrive at
\begin{eqnarray}
&&\left( \int_r^T|f(t,s)|^pds\right) ^{\frac 2p}  \nonumber \\
&\leq &\left( \frac{2-p}p\right) ^{\frac{2-p}p}\left( \frac 1\beta \right) ^{%
\frac{2-p}p}e^{-\beta A^{*}(r)}\int_r^Te^{\beta A^{*}(s)}\frac{|f(t,s)|^2}{%
\alpha ^2(s)}ds.
\end{eqnarray}
Then we have
\begin{theorem}
Let (H1) hold, $\psi(\cdot)\in L^{2,\beta}_{\mathcal{F}_{T}}[0,T],$
$\alpha ^2(s)\geq1,$ and $L(t,s)$ satisfies:
\begin{equation}
\sup\limits_{t\in [0,T]}\left( \int_t^TL^q(t,s)ds\right) ^{\frac
2q}<\infty ,
\end{equation}
where $q$ is a constant and $q>2,$ $r_i(s)$ are determined functions, $%
A^{*}(s)$ is bounded (or continuous), then BSVIE (2) admits a unique
M-solution in $\mathcal{H}^2[0,T]$.
\end{theorem}
\begin{proof}Thanks to (24), we can replace $A(s)$ in Lemma 3.1 with $A^{*}(s)$ and then (20) becomes:
\begin{eqnarray}
&&\ E\int_0^Te^{\beta A^{*}(t)}|Y(t)|^2dt+E\int_0^T\int_0^Te^{\beta
A^{*}(s)}|Z(t,s)|^2dsdt  \nonumber \\
\ &\leq &2E\int_0^Te^{\beta
A^{*}(t)}|Y(t)|^2dt+E\int_0^T\int_t^Te^{\beta
A^{*}(s)}|Z(t,s)|^2dsdt  \nonumber \\
\ &\leq &CE\int_0^Te^{\beta A^{*}(t)}|\psi
(t)|^2dt+CE\int_0^Te^{\beta
A^{*}(t)}\left| \int_t^Tf(t,s)ds\right| ^2dt  \nonumber \\
&&+CE\int_0^T|\psi (t)|^2\left( \int_t^Te^{\beta A^{*}(u)}d\beta
A^{*}(u)\right) dt  \nonumber \\
&&+\ CE\int_0^T\int_t^Te^{\beta A^{*}(s)}\left|
\int_s^Tf(t,u)du\right|
^2d\beta A^{*}(s)dt  \nonumber \\
\ &\leq &CEe^{\beta A^{*}(T)}\int_0^T|\psi (t)|^2dt+\frac C\beta
E\int_0^T\int_t^Te^{\beta A^{*}(s)}\frac{|f(t,s)|^2}{\alpha
^2(s)}dsdt.
\end{eqnarray}
The former part is the same as the corresponding part in Theorem
3.1, thus we only state the rest. Note that $Y(\cdot
),\overline{Y}(\cdot ),Z(\cdot ,\cdot ),\overline{Z}(\cdot ,\cdot )$
have the same meaning as above. It follows that,
\begin{eqnarray*}
&&E\int_0^Te^{\beta A^{*}(t)}|Y(t)-\overline{Y}(t)|^2dt+E\int_0^T\int_0^Te^{%
\beta A^{*}(s)}|Z(t,s)-\overline{Z}(t,s)|^2dsdt \\
\ &\leq &CE\int_0^Te^{\beta A^{*}(t)}\left| \int_t^TH(t,s)ds\right|
^2dt+CE\int_0^T\int_t^Te^{\beta A^{*}(s)}\left|
\int_s^TH(t,u)du\right|
^2d\beta A^{*}(s)dt \\
\ &\leq &CE\int_0^Te^{\beta A^{*}(t)}\left| \int_t^TG(t,s)ds\right|
^2dt+CE\int_0^T\int_t^Te^{\beta A^{*}(s)}\left|
\int_s^TG(t,u)du\right|
^2d\beta A^{*}(s)dt \\
&\leq &CE\int_0^Te^{\beta A^{*}(t)}\left( \int_t^TL^q(t,s)ds\right)
^{\frac 2q}\left( \int_t^TU^p(t,s)ds\right) ^{\frac 2p}dt \\
&&+CE\int_0^T\int_t^Te^{\beta A^{*}(s)}\left(
\int_s^TL^q(t,u)du\right) ^{\frac 2q}\left(
\int_s^TU^p(t,u)du\right) ^{\frac 2p}d\beta A^{*}(s)dt \\
\ &\leq &C\left( \frac 1\beta \right) ^{\frac{2-p}p}E\int_0^T\int_t^Te^{%
\beta A^{*}(s)}|y(s)-\overline{y}(s)|^2dsdt
\end{eqnarray*}
\begin{eqnarray*}
&&+C\left( \frac 1\beta \right)
^{\frac{2-p}p}E\int_0^T\int_t^Te^{\beta
A^{*}(s)}|z(t,s)-\overline{z}(t,s)|^2dsdt \\
&&+C\left( \frac 1\beta \right)
^{\frac{2-p}p}E\int_0^T\int_t^Te^{\beta
A^{*}(s)}|z(s,t)-\overline{z}(s,t)|^2dsdt \\
\ &\leq &C\left( \frac 1\beta \right)
^{\frac{2-p}p}E\int_0^Te^{\beta
A^{*}(s)}|y(s)-\overline{y}(s)|^2ds \\
&&+C\left( \frac 1\beta \right)
^{\frac{2-p}p}E\int_0^T\int_t^Te^{\beta
A^{*}(s)}|z(t,s)-\overline{z}(t,s)|^2dsdt,
\end{eqnarray*}
where
\begin{eqnarray*}
H(t,s) &=&g(t,s,y(s),z(t,s),z(s,t))-g(t,s,\overline{y}(s),\overline{z}(t,s),%
\overline{z}(s,t)), \\
G(t,s) &=&L(t,s)(r_1(s)|y(s)-\overline{y}(s)| \\
&&\ +r_2(s)|z(t,s)-\overline{z}(t,s)|+r_3(s)|z(s,t)-\overline{z}(s,t)|), \\
U(t,s) &=&r_1(s)|y(s)-\overline{y}(s)| \\
&&+r_2(s)|z(t,s)-\overline{z}(t,s)|+r_3(s)|z(s,t)-\overline{z}(s,t)|.
\end{eqnarray*}
Choosing a sufficient large number $\beta $ so that $C\left( \frac
1\beta
\right) ^{\frac{2-p}p}<1.$ Then the mapping $\Theta$ is contracted from $\mathcal{H}%
 ^{2,\beta}[0,T]$ onto itself. Because  of the assumption on
$A^{*}(s),$ we know that BSVIE (2) admits a unique M-solution in
$\mathcal{H}^2[0,T]$.
\end{proof}
\begin{remark}
In \cite{Y2}, the assumption on $L(t,s)$ is $\sup\limits_{t\in
[0,T]}\displaystyle\int_t^TL^{2+\epsilon }(t,s)ds<\infty ,$ where
$\epsilon $ is a positive constant. When $r_i(s)$ $(i=1,2,3)$ are
constants, then from (26) we know that the assumptions on Lipschitz
coefficients are the same as the one in \cite{Y2}. On the other
hand, here we assume that $A(t)$ is bounded, so if $\psi(\cdot)\in
L_{\mathcal {F}_{T}}^{2}[0,T]$ and
$E\displaystyle\int_{0}^{T}\displaystyle\int_{t}^{T}g_{0}(t,s)dsdt<\infty,$
we can get $\psi(\cdot)\in L_{\mathcal {F}_{T}}^{2,\beta}[0,T]$ and
$E\displaystyle\int_{0}^{T}\displaystyle\int_{t}^{T}e^{A(s)}g_{0}(t,s)dsdt<\infty,$
then by Theorem 3.2 we can also get the existence and uniqueness of
M-solution in $\mathcal {H}^{2}[0,T]$ which is one of the main
results in \cite{Y2}.
\end{remark}
\begin{remark}
We can use same argument as above to give the stability estimate,
for any $S\in[0,T],$
\begin{eqnarray}
&& E\int_S^T|Y(t)-\overline{Y}(t)|^2dt+\int_S^T\int_S^T|Z(t,s)-\overline{Z}%
(t,s)|^2dsdt  \nonumber \\
 &\leq& CE\int_S^T|\psi (t)-\overline{\psi
}(t)|^2dt+
CE\int_S^T\left( \int_t^T|g-\overline{g}%
|ds\right) ^2dt,
\end{eqnarray}
where $g=g(t,s,Y(s),Z(t,s),Z(s,t))$ and
$\overline{g}=\overline{g}(t,s,Y(s),Z(t,s),Z(s,t)),$ $\overline{\psi }%
(\cdot )\in L_{\mathcal{F}_T}^2[0,T]$ and $(\overline{Y}(\cdot ),\overline{Z}%
(\cdot ,\cdot ))\in \mathcal{H}^2[0,T]$ be the adapted M-solution of
(2) with $g$ and $\psi (\cdot )$ replaced by $\overline{g}$ and $\overline{%
\psi }(\cdot )$, respectively.
\end{remark}
 Similarly we can obtain the existence and uniqueness of the adapted solution for BSVIE (3).
\begin{theorem}
Let (H1) hold, $\psi(\cdot)\in L^2_{\mathcal{F}_{T}}[0,T],$ $\alpha
^2(s)\geq1,$ and $L(t,s)$ satisfies:
\[
\sup\limits_{t\in [0,T]}\left( \int_t^TL^q(t,s)ds\right) ^{\frac
2q}<\infty ,
\]
where $q$ is a constant and $q>2,$ $r_i(s)$ are two adapted processed, $%
A^{*}(s) $ is bounded, then BSVIE (3) admits a unique adapted
solution in $\mathcal{H}^2[0,T]$.
\end{theorem}
\begin{proof}
 We can get the result by Lemma 3.1 and the fixed point
 theorem, so we omit it.
 \end{proof}
\begin{remark}
Here we let the Lipschitz coefficient be stochastic because we do
not need to consider the value of $Z(t,s)$ $(0\leq s\leq t\leq T)$
for the adapted solution of BSVIE (3).
\end{remark}
\subsection{The non-Lipschitz case}
In this subsection, we will consider the unique existence of adapted
solution of BSVIE (3) and M-solution of BSVIE (2) under
non-Lipschitz condition. We assume that

(H2) For all $y,$ $\overline{y}\in R^m,$ $z,$ $\overline{z},$
$\zeta,$ $\overline{\zeta}\in R^{m\times
d}, $ and $(t,s)\in \Delta ^c$%
\begin{eqnarray*}
&&\ |g(t,s,y,z,\zeta)-g(t,s,\overline{y},\overline{z},\overline{\zeta})| \\
\ &\leq &L(t,s)(r_1(s)(\rho (|y-\overline{y}|^2))^{\frac 12}+r_2(s)|z-%
\overline{z}|+r_3(s)|\zeta-\overline{\zeta}|),
\end{eqnarray*}
where $\rho $ is an increasing concave function from $R_{+}$ to
$R_{+}$ such that $\rho (0)=0,$ and $\int_{0_{+}}\frac{du}{\rho
(u)}=\infty.$ $L(t,s)$ is a deterministic non-negative function.

Since $\rho$ is concave and $\rho(0)=0,$ one can find a pair of
positive constants $a$ and $b$ such that $\rho(u)\leq a+bu$, for all
$u\geq 0.$ Next we use the argument in Lemma 3.1 to give another
estimate which plays a critical role in the next. The following form
of BSVIE
\begin{equation}
Y(t)=\xi +\int_t^Tg(t,s,Y(s),Z(t,s))ds-\int_t^TZ(t,s)dW(s),
\end{equation}
under non-Lipschitz coefficient was considered in \cite{WZ}. The
author also gave a critical estimate by using the It\^o formula to
$e^{\beta t}|Y(t)|^2$ to give an estimate for
$e^{\beta t}|Y(t)|^2+E^{\mathcal{F}_{t}}\int^T_te^{\beta s}|Z(t,s)|^2ds$ where $%
(Y(\cdot ),Z(\cdot ,\cdot ))$ is the adapted solution of equation
(29). Now we give another estimate for $Ee^{\beta
t}|Y(t)|^2+E\int^T_te^{\beta s}|Z(t,s)|^2ds$ by the same method in
Lemma 3.1 without involving It\^o formula. We have:
\begin{lemma}
Let $\psi (\cdot )\in L_{\mathcal{F}_{T} }^{2,\beta}[0,T],$ the
assumptions on $f$ is the same as in Lemma 3.1,
$(Y(\cdot),Z(\cdot,\cdot))$ is the adapted solution of (4), then for
almost every $t\in [0,T],$ we have the following estimate:
\begin{eqnarray*}
&&Ee^{\beta A(t)}|Y(t)|^2+E\int_t^Te^{\beta A(s)}|Z(t,s)|^2ds \\
&\leq& Ee^{\beta A(T)}|\psi (t)|^2+\frac C\beta E\int_t^Te^{\beta A(s)}\frac{%
|f(t,s)|^2}{\alpha ^2(s)}ds.
\end{eqnarray*}
\end{lemma}
\begin{proof}It follows from (7) that
\begin{eqnarray*}
Y(t) &=&\lambda (t,t)=E^{\mathcal{F}_t}\left\{ \psi
(t)+\int_t^Tf(t,s)ds\right\}.
\end{eqnarray*}
Then
\begin{eqnarray*}
e^{\beta A(t)}|Y(t)|^2 &\leq &2e^{\beta A(t)}E(\left. |\psi
(t)|^2\right| \mathcal{F}_t)+2e^{\beta A(t)}E\left( \left. \left|
\int_t^Tf(t,s)ds\right| ^2\right| \mathcal{F}_t\right) ,
\end{eqnarray*}
thus
\begin{eqnarray*}
Ee^{\beta A(t)}|Y(t)|^2 &\leq &2E(e^{\beta A(t)}|\psi (t)|^2)
+2E\left( e^{\beta A(t)}\left| \int_t^Tf(t,s)ds\right| ^2\right).
\end{eqnarray*}
By taking $r=t$ in (11), we claim that
\begin{eqnarray*}
&&E\int_t^Te^{\beta A(s)}|Z(t,s)|^2ds \\
&=&E\int_t^Te^{\beta A(s)}\left( \int_s^T|Z(t,u)|^2du\right) d\beta
A(s)
 +Ee^{\beta A(t)}\int_t^T|Z(t,u)|^2du.
\end{eqnarray*}
In the light of the proof in Lemma 3.1, it can be easily checked
that,
\begin{eqnarray*}
&&E\int_t^Te^{\beta A(s)}|Z(t,s)|^2ds \\
&\leq &CEe^{\beta A(T)}|\psi (t)|^2+CE\left( e^{\beta A(t)}\left|
\int_t^Tf(t,s)ds\right| ^2\right) \\
&&+CE\left(\int_t^Te^{\beta A(s)}\left| \int_s^Tf(t,u)du\right|
^2d\beta A(s)\right).
\end{eqnarray*}
Then we obtain that
\begin{eqnarray*}
&&Ee^{\beta A(t)}|Y(t)|^2+E\int_t^Te^{\beta A(s)}|Z(t,s)|^2ds \\
&\leq &CEe^{\beta A(T)}|\psi (t)|^2+CE\left( e^{\beta A(t)}\left|
\int_t^Tf(t,s)ds\right| ^2\right) \\
&&+CE\left(\int_t^Te^{\beta A(s)}\left| \int_s^Tf(t,u)du\right| ^2d\beta A(s)\right) \\
&\leq &CEe^{\beta A(T)}|\psi (t)|^2+\frac C\beta E\int_t^Te^{\beta A(s)}\frac{%
|f(t,s)|^2}{\alpha ^2(s)}ds.
\end{eqnarray*}
The conclusion thus follows.
\end{proof}
\begin{remark}
When $r_i(s)(i=1,2)$ are constants, the above estimate becomes:
\begin{eqnarray}
\ Ee^{\beta t}|Y(t)|^2+E\int_t^Te^{\beta s}|Z(t,s)|^2ds \nonumber \\
\ \leq CEe^{\beta T}|\psi (t)|^2+\frac C\beta E\int_t^Te^{\beta
s}|f(t,s)|^2ds
\end{eqnarray}
which is similar to the one in \cite{WZ}:
\begin{eqnarray}
\ e^{\beta t}|Y(t)|^2+E^{\mathcal{F}_t}\int_t^Te^{\beta s}|Z(t,s)|^2ds  \nonumber \\
\ \leq e^{\beta T}E^{\mathcal{F}_t}|\xi |^2+\frac 1\beta
E\int_t^Te^{\beta s}|f(t,s)|^2ds.
\end{eqnarray}
Though the estimate (31) is stronger than (30), (30) still can
guarantee that equation (3) admits a unique adapted solution under
non-Lipschitz coefficients.
\end{remark}
We have
\begin{theorem}
Let (H2) hold, $g$ is independent of $Z(s,t)$, $\psi(\cdot)\in
L^2_{\mathcal{F}_{T}}[0,T]$, $r_i(s)$ are adapted processes, $A(t)$
is bounded, $L(t,s)$ satisfies:
\[
\sup\limits_{t\in [0,T]}\left( \int_t^TL^q(t,s)ds\right) ^{\frac
2q}<\infty ,
\]
where $q>2$ is a constant, then (?) admits a unique adapted solution
in $\mathcal{H}^2_t[0,T]$.
\end{theorem}
\begin{proof}The proof can be obtained by combining the estimate in Lemma
3.1 and the proof in \cite{WZ}, so we omit it.
\end{proof}

When $r_i(s)$ is a constant and $\psi (\cdot )=\xi,$ $g$ is
independent of $Z(s,t),$ then we get the result in \cite{WZ}:
\begin{corollary}
Let (H2) hold, $r_i(s)=1$, $L(t,s)=k,$ $k$ is a constant, then (3)
admits a unique adapted solution in $\mathcal{H}^2_t[0,T]$.
\end{corollary}

Obviously we can also get the unique existence of M-solution of (3)
as above. However, as to the general form of (2), there is an
expression $E\int_t^T|Z(s,t)|^2ds$ which is hard to estimate
directly, see \cite{Y2}, so we have to adopt new method to deal with
it. We need to prepare some results in order to be able to derive
this claim.
\begin{lemma}
 For any $t\in [0,T],$
$f(s):[t,T]\rightarrow R^{+},$ $c(x):R\rightarrow R$ is a concave
function, $\int_0^Tf(s)ds<\infty.$ Then we have
\[
\frac 1{T-t}\int_t^Tc(f(s))ds\leq c\left( \frac
1{T-t}\int_t^Tf(s)ds\right) .
\]
\end{lemma}
\begin{proof}
Obviously $-c(x)$ is a convex function, for fixed $x\in R,$ $\forall $ $%
y_1>x,$ $y_2<x,$ we have (see \cite{RY}).
\[
\frac{-c(y_1)+c(x)}{y_1-x}\geq -c_{+}^{^{\prime }}(x)\geq
-c_{-}^{^{\prime }}(x)\geq \frac{-c(y_2)+c(x)}{y_2-x},
\]
thus there exists a $k\in [-c_{-}^{^{\prime }}(x),-c_{+}^{^{\prime
}}(x)]$
so that $\forall y\in R,$ $-c(y)\geq -c(x)+k\cdot(y-x),$ i.e., $%
c(y)\leq c(x)-k(y-x).$ For any fixed $t\in [0,T],$ $s\in [t,T],$
\[
x=\frac 1{T-t}\int_t^Tf(s)ds,\text{ }y=f(s),
\]
then
\[
c(f(s))\leq c\left( \frac 1{T-t}\int_t^Tf(s)ds\right) -k\cdot\left(
f(s)-\frac 1{T-t}\int_t^Tf(s)ds\right) ,
\]
thus we get the conclusion above.
\end{proof}

We are now ready to establish the last result of this paper.
\begin{theorem}
Let (H2) hold, $\psi(\cdot)\in L^2_{\mathcal{F}_{T}}[0,T]$, $r_i(s)$
are deterministic function, $A(t)$ is bounded, $L(t,s)$ satisfies:
\[
\sup\limits_{t\in [0,T]}\left( \int_t^TL^q(t,s)ds\right) ^{\frac
2q}<\infty ,
\]
where $q>2$ is a constant, then (2) admits a unique M-solution in
$\mathcal{H}^2[0,T]$.
\end{theorem}
\begin{proof}
Uniqueness: Let $(Y_i,Z_i)\in \mathcal{H}^2[0,T]$ $(i=1,2)$ be any
two M-solutions. By defining
\[
\widehat{Y}(t)=Y_1(t)-Y_2(t);\text{ }\widehat{Z}(t,s)=Z_1(t,s)-Z_2(t,s),%
\text{ }t,s\in [0,T],
\]
we arrive that
\begin{eqnarray*}
&&\widehat{Y}(t)+\int_t^T\widehat{Z}(t,s)dW(s) \\
&=&%
\int_t^T[g(t,s,Y_1(s),Z_1(t,s),Z_1(s,t))-g(t,s,Y_2(s),Z_2(t,s),Z_2(s,t))]ds.
\end{eqnarray*}
Note that $\widehat{Y}(T)=0.$ For arbitrary $u\in [0,T),$ we obtain
the following results in the same way as Theorem 3.2,
\begin{eqnarray*}
&&E\int_u^Te^{\beta
A^{*}(t)}|\widehat{Y}(t)|^2dt+E\int_u^T\int_t^Te^{\beta
A^{*}(s)}|\widehat{Z}(t,s)|^2dsdt \\
&\leq &C\left( \frac 1\beta \right)
^{\frac{2-p}p}E\int_u^T\int_t^Te^{\beta
A^{*}(s)}\rho (|\widehat{Y}(s)|^2)dsdt \\
&&+C\left( \frac 1\beta \right)
^{\frac{2-p}p}E\int_u^T\int_t^Te^{\beta
A^{*}(s)}|\widehat{Z}(t,s)|^2dsdt+C\left( \frac 1\beta \right) ^{\frac{2-p}p}E\int_u^Te^{\beta A^{*}(t)}|%
\widehat{Y}(t)|^2dt.
\end{eqnarray*}
By choosing a suitable $\beta $, we deduce the following
\[
E\int_u^Te^{\beta A^{*}(t)}|\widehat{Y}(t)|^2dt\leq
CE\int_u^T\int_t^Te^{\beta A^{*}(s)}\rho (|\widehat{Y}(s)|^2)dsdt,
\]
consequently,
\begin{eqnarray*}
\frac 1{T-u}E\int_u^T|\widehat{Y}(t)|^2dt &\leq &CE\int_u^T\frac
1{T-t}\int_t^T\rho (|\widehat{Y}(s)|^2)dsdt \\
&\leq &C\int_u^T\rho \left( \frac
1{T-t}\int_t^TE|\widehat{Y}(s)|^2ds\right) dt.
\end{eqnarray*}
Due to Bihari's inequality (see \cite{B}) we get that $\frac
1{T-u}E\int_u^T|\widehat{Y}(t)|^2dt=0,$
$u\in [0,T),$ thus $\widehat{Y}(t)=0\ $as well as $\widehat{Z}(t,s)=0,$ $%
t,s\in [0,T].$ a.e.

Existence: Let $Y_0(t)\equiv 0,$ and define recursively $(Y_n,Z_n)$
by the following equations with the help of Theorem 3.2
\begin{eqnarray}
Y_n(t)=\psi
(t)+\int_t^Tg(t,s,Y_{n-1}(s),Z_n(t,s),Z_n(s,t))ds-\int_t^TZ_n(t,s)dW(s).%
\end{eqnarray}
By setting
\[
\widehat{Y}_{n,k}(t)=Y_n(t)-Y_k(t);\text{ }\widehat{Z}%
_{n,k}(t,s)=Z_n(t,s)-Z_k(t,s),\text{ }t,s\in [0,T],
\]
and choosing a suitable $\beta ,$ we claim that
\begin{eqnarray*}
&&E\int_u^Te^{\beta A^{*}(t)}|\widehat{Y}_{n,k}(t)|^2dt+E\int_u^T\int_t^Te^{%
\beta A^{*}(s)}|\widehat{Z}_{n,k}(t,s)|^2dsdt \\
&\leq &CE\int_u^T\int_t^Te^{\beta A^{*}(s)}\rho (|\widehat{Y}%
_{n-1,k-1}(s)|^2)dsdt,
\end{eqnarray*}
then
\[
\frac 1{T-u}E\int_u^T|\widehat{Y}_{n,k}(t)|^2dt\leq C\int_u^T\rho
\left( \frac 1{T-t}\int_t^TE|\widehat{Y}_{n-1,k-1}(s)|^2ds\right)
dt.
\]
Set $Q(u)=\limsup\limits_{n,k\rightarrow \infty }E\int_u^T|\widehat{Y}%
_{n,k}(t)|^2dt,$ it is easy to show that $S(u)=\sup\limits_{n\geq
0}E\int_u^T|Y_n(t)|^2dt$ is bounded. In fact, using the similar
trick as in Theorem 3.2, we obtain that
\begin{eqnarray*}
&&E\int_u^Te^{\beta A^{*}(t)}|Y_n(t)|^2dt+E\int_u^T\int_t^Te^{\beta
A^{*}(s)}|Z_n(t,s)|^2dsdt \\
&\leq &CE\int_u^Te^{\beta A^{*}(t)}|\psi (t)|^2dt+C\left( \frac
1\beta
\right) ^{\frac{2-p}p}E\int_u^T\int_t^Te^{\beta A^{*}(s)}|g_0(t,s)|^2dsdt \\
&&+C\left( \frac 1\beta \right)
^{\frac{2-p}p}E\int_u^T\int_t^Te^{\beta
A^{*}(s)}(a+b|Y_{n-1}(s)|^2)dsdt \\
&&+C\left( \frac 1\beta \right)
^{\frac{2-p}p}E\int_u^T\int_t^Te^{\beta A^{*}(s)}|Z_n(t,s)|^2dsdt
+C\left( \frac 1\beta \right) ^{\frac{2-p}p}E\int_u^Te^{\beta
A^{*}(t)}|Y_n(t)|^2dt,
\end{eqnarray*}
thus by choosing $%
\beta $ we have
\begin{eqnarray*}
E\int_u^T|Y_n(t)|^2dt &\leq &C+CE\int_u^T|\psi (t)|^2dt \\
&&+CE\int_u^T\int_t^T|g_0(t,s)|^2dsdt+E\int_u^T\int_t^T|Y_{n-1}(s)|^2dsdt.
\end{eqnarray*}
In view of Gronwall's inequality we obtain that $S(u)$is bounded.
Then by Fatou's lemma, Bihari's inequality and noting that $\rho $
is increasing, we deduce that for almost $u\in [0,T],$ $Q(u)=0,$ and
it follows that
\[
\lim\limits_{n,k\rightarrow \infty }E\int_0^T|Y_n(t)-Y_k(t)|^2dt=0,
\]
hence there is a $Y$ such that
\[
\lim\limits_{n\rightarrow \infty }E\int_0^T|Y_n(t)-Y(t)|^2dt=0.
\]
Similarly there is a $Z$ such that
\begin{eqnarray*}
\lim\limits_{n\rightarrow \infty
}E\int_0^T\int_t^T|Z_n(t,s)-Z(t,s)|^2dsdt
&=&0, \\
\lim\limits_{n\rightarrow \infty
}E\int_0^T\int_0^t|Z_n(t,s)-Z(t,s)|^2dsdt &\leq
&\lim\limits_{n\rightarrow \infty }E\int_0^T|Y_n(t)-Y(t)|^2dt=0.
\end{eqnarray*}
By taking the limits for BSVIE (32), one can finds that $(Y,Z)$ is a
M-solution of BSVIE (2).
\end{proof}

At last we want to give a simple example to show the unique
existence of adapted solution (or M-solution) under non-Lipschitz
condition. As shown in \cite{M1} or \cite{M2}, the following two
functions satisfy the assumption of $\rho$ in (H2),
\[
\rho _1(x)=\left\{
\begin{array}{cc}
x\ln (x^{-1}), & x\in [0,\delta ], \\
 \delta \ln (\delta ^{-1})+\dot{\rho
_1}(\delta -)(x-\delta ), & x>\delta ,
\end{array}
\right.
\]
\[
\rho _2(x)=\left\{
\begin{array}{cc}
x\ln (x^{-1})\ln \ln (x^{-1}), & x\in [0,\delta ], \\
\delta \ln (\delta ^{-1})\ln \ln (\delta ^{-1})+\dot{\rho _2}(\delta
-)(x-\delta ), & x>\delta ,
\end{array}
\right.
\]
with $\delta \in (0,1)$ being sufficiently small. However, the
explicit form of $\rho_i$ is not easy to get, so now we will give
another example to avoid this problem.

Let us consider the following equation
\begin{eqnarray}
Y(t)=\psi(t)+\int_t^TL(t,s)[f(|Y(s)|)+|Z(t,s)|+|Z(s,t)|]ds-\int_t^TZ(t,s)dW(s),
\end{eqnarray}
where $f:R\rightarrow [0,\infty)$ is defined by
\[
f(x)=\left\{
\begin{array}{cc}
0 & x=0, \\
|x|\left[ \ln (1+|x|^{-1})\right] ^{\frac 12} & 0<|x|<\delta , \\
\delta \left[ \ln (1+|\delta |^{-1})\right] ^{\frac 12} & |x|\geq
\delta ,
\end{array}
\right.
\]
$L(t,s)$ satisfies
$\sup\limits_{t\in[0,T]}\int_t^TL^{q}(t,s)ds<\infty,$ where $q>2$ is
a constant. It can be shown that $|f(y)-f(\overline{y})|\leq
\rho(|y-\overline{y}|^2)^{\frac 1 2}$ where $\rho$ can be defined by
\[
\rho(x)=\left\{
\begin{array}{cc}
0 & x=0, \\
x\ln (1+x^{-1}) & 0<x<1, \\
\ln 2 & x\geq 1,
\end{array}
\right.
\]
We refer the reader to \cite{C} for the proof. Then by Theorem 3.5,
we deduce that BSVIE (33) admits a unique M-solution. Note that we
can give the example for adapted solution in a similar way.

\end{document}